\documentclass{article}%
\usepackage{amsmath}%
\usepackage{amsfonts}%
\usepackage{amssymb}%
\usepackage{graphicx}
\providecommand{\U}[1]{\protect\rule{.1in}{.1in}}
\newtheorem{theorem}{Theorem}

\newtheorem{example}[theorem]{Example}

\begin{document}

\title{Existence of reciprocal matrices with specified orders for the  right and
inverse left Perron eigenvectors}
\author{Susana Furtado\thanks{Corresponding author. Email: sbf@fep.up.pt.
orcid.org/0000-0003-0395-5972. The work of this author is funded by FCT -
Funda\c{c}\~{a}o para a Ci\^{e}ncia e a Tecnologia, I.P., through national
funds, under the project UID/04561/2025,
https://doi.org/10.54499/UID/04561/2025.}\\CEMS.UL and Faculdade de Economia \\Universidade do Porto\\Rua Dr. Roberto Frias\\4200-464 Porto, Portugal
\and Charles R. Johnson \thanks{Email: crjmatrix@gmail.com. }\\225 West Tazewells Way\\Williamsburg, VA 23185, United States}
\maketitle

\begin{abstract}
Here we give a procedure to construct a reciprocal matrix for which the right
and entrywise inverse left Perron eigenvectors have any pair of given orders.
An explicit example when the matrix is of size $4$ is presented. In
particular, it gives an afirmative answer to the question posed in a recent
manuscript by Boz\'{o}ki and Csat\'{o} (2026) about the existence of a
reciprocal matrix of size $4$ such that the right and entrywise inverse left
Perron eigenvectors have reverse orders.

\end{abstract}

An $n$-by-$n$ \emph{pairwise comparison} matrix, or \emph{reciprocal} matrix,
$A=[a_{ij}]$ is an entrywise positive matrix in which $a_{ji}=\frac{1}{a_{ij}%
},$ $1\leq i,j\leq n$. The diagonal entries are all $1$, and the off-diagonal
entries represent pairwise ratio comparisons among $n$ alternatives. Let
$\mathcal{PC}_{n}$ denote the set of all such matrices. If $A\in
\mathcal{PC}_{n}$ satisfies, in addition, $a_{ij}a_{jk}=a_{ik}$, $1\leq
i,j,k\leq n$, then $A$ is further said to be \emph{consistent}.

\bigskip

One of the most used weight vectors obtained from a reciprocal matrix is the
right Perron eigenvector proposed in \cite{saaty1977}. Later, in \cite{John}
it was observed that the entrywise inverse of the left Perron eigenvector
could be used as well. However, an example was given in which the two vectors
have different orders. In \cite{Csato} an example of a  $5$-by-$5$  reciprocal
matrix for which the right and entrywise inverse left Perron eigenvectors have
reverse orders  is given. In the recent manuscript \cite{Boz} it is claimed
that no such example is known for $4$-by-$4$ reciprocal matrices. 

\bigskip

The right and entrywise inverse left Perron eigenvectors of a matrix in
$\mathcal{PC}_{3}$ are known to coincide \cite{Boz}. Here we give a way to
construct a matrix in $\mathcal{PC}_{n}$, $n\geq4$, for which the right and
entrywise inverse left Perron eigenvectors have any pair of given orders.
Moreover, we show that such a matrix can be arbitrarily close to a consistent
matrix. We illustrate our method in the $4$-by-$4$ case. In particular, we
show that the two above mentioned vectors may have reverse orders, answering a
question posed in \cite{Boz}. 

\bigskip

Let $A\in\mathcal{PC}_{n}$ be inconsistent. We construct a reciprocal matrix,
positive diagonal similar to $A,$ such that its right and (entrywise inverse)
left Perron eigenvectors have any desired orders. The only restriction we
impose on $A$ is that, if $A^{\prime}$ is the matrix diagonal similar to $A$
with constant row sums (see \cite{FJ5} for more details), then the left Perron
eigenvector of $A^{\prime}$ has distinct entries. This is the generic
situation. In our construction, we may assume that one of the orders is fixed
(for example increasing), as some permutation similarity on the matrix changes
that order to any desired one. 

\bigskip

Let $A\in\mathcal{PC}_{n}$ and let $v$ be its right Perron eigenvector. Let
$S=\operatorname*{diag}(v)$ be the positive diagonal matrix with diagonal $v$.
Then, the right Perron eigenvector of $A^{\prime}:=S^{-1}AS$ is $S^{-1}%
v=\mathbf{e}_{n}$, the vector of ones. This means that $A^{\prime}$ has
constant row sums. Suppose that the left Perron eigenvector $u$ of $A^{\prime
}$ has distinct entries. Then, after a sufficiently small perturbation of
$A^{\prime}$ via a positive diagonal similarity, $DA^{\prime}D^{-1}$, the new
left Perron eigenvector has distinct entries and the same order as $u$. Let
$w$ be a positive $n$-vector with any desired order and with entries
arbitrarily close to $1$. Then, for $D=\operatorname*{diag}(w)$, the right
Perron eigenvector of $B:=DA^{\prime}D^{-1}$ is $D\mathbf{e}_{n}=w$. Hence,
the right Perron eigenvector of $B$ has the order of $w$, while the left
Perron eigenvector has the order of $u$. If matrix $A$ is close to
consistency, then so is $B$. Notice that $A$ and $B$ are positive diagonal
similar. So they have the same Perron root $\lambda$ and, therefore, the same
inconsistency index $CI=\frac{\lambda-n}{n-1}$ \cite{saaty1977}.

\bigskip

Next we give an example when $n=4.$

\begin{example}
Let
\[
A=\left[
\begin{array}
[c]{cccc}%
1 & 0.6801 & 0.1181 & 0.5869\\
1.4704 & 1 & 0.2605 & 1.7259\\
8.4660 & 3.8385 & 1 & 1.6563\\
1.7038 & 0.5794 & 0.6038 & 1
\end{array}
\right]  \in\mathcal{PC}_{4}.
\]
The Perron eigenvector of $A$ is $v=\left[
\begin{array}
[c]{cccc}%
0.5 & 1 & 3 & 1
\end{array}
\right]  ^{T}.$ Let $S=\operatorname*{diag}(v).$ Then
\[
A^{\prime}=S^{-1}AS=\left[
\begin{array}
[c]{cccc}%
1 & 1.3602 & 0.7087 & 1.1738\\
0.7352 & 1 & 0.781\,53 & 1.\,\allowbreak725\,9\\
1.4110 & 1.2795 & 1 & 0.5521\\
0.8519 & 0.5794 & 1.8113 & 1
\end{array}
\right]  .
\]
The right Perron eigenvector of $A^{\prime}$ is $\mathbf{e}_{4}$. The left
Perron eigenvector of $A^{\prime}$, that is, the right Perron eigenvector of
$A^{\prime T}$, is $\ $%
\[
u=\left[
\begin{array}
[c]{cccc}%
0.9062 & 0.9475 & 0.9850\, & 1
\end{array}
\right]  ^{T}.
\]
Vector $u$ has strictly increasing entries. Thus, a small perturbation of
$A^{\prime}$ via a positive diagonal similarity leaves the left Perron
eigenvector strictly increasing. Let $w$ be a positive vector with entries
arbitrarily close to $1$ and with any given order. Let $D=\operatorname*{diag}%
(w)$. Then, the right Perron eigenvector of $B:=DA^{\prime}D^{-1}$ is $w$ and
the left Perron eigenvector of $B$ has the same order as $u$. 

For example, if  $w=(0.97,0.98,0.99,1)$, the left Perron eigenvector of $B$
is
\[
\left[
\begin{array}
[c]{cccc}%
0.9343 & 0.9669 & 0.995\, & 1
\end{array}
\right]  ^{T}.
\]
In this case, the right and entrywise inverse left Perron eigenvectors of $B$
have reverse orders.

If $w=(1,0.99,0.98,0.97)$, the left Perron eigenvector of $B$ is
\[
\left[
\begin{array}
[c]{cccc}%
0.8790\, & 0.9283 & 0.9750\, & 1
\end{array}
\right]  ^{T},
\]
and the right and entrywise inverse left Perron eigenvectors of $B$ have the
same order.
\end{example}

\end{document}